\documentclass[12pt,draft]{amsart}
\usepackage{amsmath,amsthm,latexsym,amscd,amsbsy,amssymb,pb-diagram}
 \relax

\chardef\bslash=`\\ 

\makeatletter
\def\verbatim{\interlinepenalty\@M \@verbatim
  \leftskip\@totalleftmargin\advance\leftskip2pc
  \frenchspacing\@vobeyspaces \@xverbatim}
\makeatother
\makeatletter
  \def\dgt@k{\dg@DX=1 \dg@DY=-4 \dg@SIZE=1}
\makeatother

\theoremstyle{plain}
\newtheorem{thm}{Theorem}[section]
\newtheorem{cor}[thm]{Corollary}
\newtheorem{lem}[thm]{Lemma}
\newtheorem{pro}[thm]{Proposition}

\newtheorem*{A}{Theorem A}
\newtheorem*{B}{Theorem B}
\newtheorem*{C}{Theorem C}
\newtheorem*{D}{Theorem D}

\theoremstyle{definition}

\newtheorem{defin}[thm]{Definition}
\newtheorem{example}[thm]{Example}

\newcounter{rmnum}

\newcommand{\grup}[1]{\pi_{#1}^{[L]}}
\newcommand{\disk}[1]{D^{#1}_{[L]}}
\newcommand{\sfr}[1]{S^{#1}_{[L]}}
\newcommand{\tdisk}[1]{\widetilde{D^{#1}_{[L]}}}
\newcommand{\ddisk}[1]{{\mathbb D} ^{#1}_{[L]}}
\newcommand{\ssfr}[1]{{\mathbb S} ^{#1}_{[L]}}
\newcommand{\cyl}[1]{C^{#1}_{[L]}}
\newcommand{\point}{s_{[L]}}
\newcommand{\pt}{s_{[L]}}
\def\ed{\text{e}-\dim}
\def\Int{{\rm Int}}

\def\fL{f^{[L]}}

\numberwithin{equation}{section}

\font\f=msbm10

\begin{document}


\title[Topological model categories generated by finite complexes]
{Topological model categories generated by finite complexes}
\author{Alex Chigogidze}
\address{Department of Mathematics and Statistics,
University of Saskatche\-wan,
McLean Hall, 106 Wiggins Road, Saskatoon, SK, S7N 5E6,
Canada}
\email{chigogid@math.usask.ca}
\thanks{The first author was partially supported by NSERC research grant.}
\author{A. Karasev}
\address{Department of Mathematics and Statistics,
University of Saskatche\-wan,
McLean Hall, 106 Wiggins Road, Saskatoon, SK, S7N 5E6,
Canada}
\email{karasev@math.usask.ca}



\begin{abstract}{Our main result states that for each finite complex $L$ the category ${\bf TOP}$ of topological spaces possesses a model category structure (in the sense of Quillen) whose weak equivalences are precisely maps which induce isomorphisms of all $[L]$-homotopy groups. The concept of $[L]$-homotopy has earlier been introduced by the first author and is based on Dranishnikov's notion of extension dimension. As a corollary we obtain an algebraic characterization of $[L]$-homotopy equivalences between $[L]$-complexes. This result extends two classical theorems of J.~H.~C.~Whitehead. One of them -- describing homotopy equivalences between CW-complexes as maps inducing isomorphisms of all homotopy groups -- is obtained by letting $L = \{ {\rm point}\}$. The other -- describing $n$-homotopy equivalences between at most $(n+1)$-dimensional CW-complexes as maps inducing isomorphisms of $k$-dimensional homotopy groups with $k \leq n$ -- by letting $L = S^{n+1}$, $n \geq 0$.}
\end{abstract}

\maketitle \markboth{A.~Chigogidze, A.~Karasev}{Topological model
categories generated by finite complexes}

\section{Introduction}\label{S:intro}
The basic concept of the model category, introduced by Quillen \cite{qui}, provides an extremely useful tool for developing axiomatic homotopy theory in very general situations (see, for instance, \cite{kan1}, \cite{kan2}, \cite{heller}, \cite{brown}). Recall that a model category structure can be introduced into a category by specifying three classes of morphisms, called {\em fibrations}, {\em cofibration} and {\em weak equivalences}, which satisfy certain axioms. By passing to the ``homotopy category" ${\mathcal Ho}({\mathcal C})$ of a model category ${\mathcal C}$ one formally inverts weak equivalences (if both domain and range are nice -- cofibrant and fibrant simultaneously). In other words, in the quotient category ${\mathcal Ho}({\mathcal C})$ weak equivalences become homotopy equivalences. This fundamental fact manifests itself in different ways in particular situations. Even though the axioms of model categories are verifiable in a wide variety \cite{dwyspa95},  \cite{bau88}, \cite{hov} of situations, they are reminiscent of well-known properties of homotopies for topological spaces. Obviously the category ${\bf TOP}$ of topological spaces itself possesses a model category structure (see, for instance, \cite{dwyspa95}).

\begin{A}
The category {\bf TOP} of topological spaces admits a model category structure where a map $f$ is a weak equivalence if $f$ is a weak homotopy equivalence, i.e. for each $k = 0,1,\dots $ and each $x \in X$ the induced map $\pi_{k}(X, x) \to \pi_{k}(Y,f(x))$ is an isomorphism.
\end{A}

It should be noted that there exists \cite{strom} a model category structure on ${\bf TOP}$ whose weak equivalences are just standard homotopy equivalences. These two structures are essentially different although the weak equivalences between CW-complexes in the corresponding homotopy categories become invertible, i.e. coincide with the ordinary homotopy equivalences. Clearly, for the latter structure this fact is true by definition (for any spaces). As for the former, the corresponding fact simply restates the following well known theorem of J.H.C.Whitehead \cite{whi1, whi2} (note that $CW$-complexes are cofibrant and fibrant in the first structure).

\begin{B}
A map between CW-complexes (or, more generally, {\em ANE}-spaces) is a homotopy equivalence if and only if it induces isomorphisms of {\em all} homotopy groups.
\end{B}

There is one more type of model category structures on ${\bf TOP}$ closely related to the ordinary homotopies. Consider, for each $n = 0,1,\dots$, Whitehead's $n$-types \cite{whi3} based on the concept of $n$-homotopy introduced in \cite{fox}. Algebraic models for $n$-types (the so called $cat^{n}$-groups) were found in \cite{loday}. Approximately at the same time $n$-homotopies (and subsequently even $n$-shapes \cite{chi88}, \cite{chi89}) begun to play a substantial role in the revitalized theory of Menger manifolds \cite{bestvina}, \cite{chibook} (see \cite{chife} for a discussion of categorical connections between $n$-homotopies and homotopies via the theories of manifolds modeled on Menger and Hilbert cubes respectively). The following theorem has been proved in \cite{donazar}.

\begin{C}
Let $n = 0,1,\dots$. The category {\bf TOP} of topological spaces admits a model category structure where a map $f \colon X \to Y$ is
a weak equivalence if it is a weak $n$-homotopy equivalence, i.e. for each $k = 0,1,\dots, n$ and each $x \in X$ the induced map $\pi_{k}(X, x) \to \pi_{k}(Y,f(x))$ is an isomorphism.
\end{C}

The corresponding homotopy category ${\mathcal Ho}_{n}({\bf TOP})$ is a model category for the above mentioned notion of $n$-type. It is essential to note that the invertibility of weak equivalences in  ${\mathcal Ho}_{n}({\bf TOP})$, in analogy with Theorem B, restates another well-known result of Whitehead \cite{whi1}.

\begin{D}
A map between at most $n$-dimensional CW-complexes (or, more generally, at most $n$-dimensional {\em LC}$^{n-1}$-spaces) is an $n$-homotopy equivalence if and only if it induces isomorphisms of the $k$-th homotopy groups for each $k \leq n-1$.
\end{D}

Below, for each finite CW-complex $L$ we consider the concept of $[L]$-homotopy introduced in \cite{chi}. We also present detailed description of $L$-homotopy groups also introduced in \cite{chi}. Our main result is as follows.
\medskip

\noindent {\bf Theorem \ref{T:main}.} {\em Let $L$ be a finite $CW$-complex. The category {\bf TOP} of topological spaces admits a model category structure where a map $f \colon X \to Y$ is
a weak equivalence if it is a weak $[L]$-homotopy equivalence, i.e. for each $n = 0,1,\dots$ and each $x \in X$ the induced map $\pi_{n}^{[L]}(X, x) \to \pi_{n}^{[L]}(Y,f(x))$ of $[L]$-homotopy groups is an isomorphism.}
\medskip

As above we conclude that weak equivalences become invertible in the corresponding homotopy category ${\mathcal Ho}_{[L]}({\bf TOP})$. This proves the earlier announced result from \cite[Theorem 2.9]{chi}.
\medskip

\noindent {\bf Corollary \ref{C:whitehead}.} {\em A map between $[L]$-complexes is an $[L]$-homotopy equivalence if and only if it induces isomorphisms of all $[L]$-homotopy groups.}
\medskip

This result also extends both theorems B and D. The first is obtained by letting $L = \{ \text{{\rm point}}\}$ and the other by assuming $L = S^{n}$.

The concept of $[L]$-homotopy differs from the notions of homotopy or $n$-homotopy (see comment right before Definition \ref{D:absoluteextensors}) and it seems very interesting to develop this theory further as well as to consider corresponding homology and cohomology theories. There are some indications that such theories could be better designed for investigation of particular geometric constructions arising within the extension theory. Conversely, particular methods used in developing the theory of extension dimension (which led to $[L]$-homotopies) could be useful in building dimension theory in particular model categories.

The paper is organized as follows. In Section \ref{S:pre} we collect necessary definitions related to model categories and extension types of complexes. Section \ref{S:lhomotopies} discusses the concept of $[L]$-homotopy, introduced earlier by the first author. In particular, we define $[L]$-homotopy groups. Finally, in Section \ref{S:model}, we prove our main result (Theorem \ref{T:main}) which provides an explicit description of a model category structure of {\bf TOP} whose weak equivalences are precisely maps inducing isomorphisms of all $[L]$-homotopy groups.

The authors are grateful to the referee whose comments and
suggestions led to a substantial improvement of the original
exposition. 


\section{Preliminaries}\label{S:pre}
In this section we present relevant definitions and results regarding model categories and extension dimension.
\subsection{Model categories}\label{SS:model} We begin with the Quillen's concept of the model category.
\begin{defin}\label{D:model}
A model category is a category ${\mathcal C}$ with
three distinguished classes of maps:
\begin{itemize}
\item[(i)]
weak equivalences,
\item[(ii)]
fibrations,
\item[(iii)]
cofibrations,
\end{itemize}
each of which is closed under compositions and contains
all identity maps. A map which is both a fibration
(respectively, cofibration) and a weak equivalence is called an
acyclic fibration (respectively, acyclic cofibration). We require
the following axioms:
\begin{itemize}
\item[(MC1)]
Finite limits and colimits exist in ${\mathcal C}$.
\item[(MC2)]
If $f$ and $g$ are maps in ${\mathcal C}$ such that $gf$ is
defined and if two of the three maps $f$, $g$, $gf$ are weak
equivalences, then so is the third.
\item[(MC3)]
If $f$ is a retract of $g$ and $g$ is a fibration, cofibration
or a weak equivalence, then so is $f$.
\item[(MC4)]
Given a commutative diagram
\[
\begin{diagram}
\node{A} \arrow{e,t}{f} \arrow{s,l}{i}  \node{X} \arrow{s,r}{p}\\
\node{B} \arrow{e,t}{g} \arrow{ne,t,..}{h} \node{Y}
\end{diagram}
\]
\noindent of unbroken arrows, a lift (the broken arrow) exists
in either of the following two situations: (i) $i$ is a
cofibration and $p$ is an acyclic fibration, or (ii) $i$ is
an acyclic cofibration and $p$ is a fibration.
\item[(MC5)]
Any map $f$ can be factored in two ways: (i) $f = pi$, where
$i$ is a cofibration and $p$ is an acyclic fibration, and (ii)
$f = qj$, where $j$ is an acyclic cofibration and $q$ is a fibration.
\end{itemize}
\end{defin}

If ${\mathcal C}$ is a model category, then it has an initial object and a terminal object (the first being the colimit and the second being the limit of the empty diagram). An object in ${\mathcal C}$ is called {\em cofibrant} if the map from the initial object to it is a cofibration. Similarly an object in ${\mathcal C}$ is called {\em fibrant} if the map from it to the terminal object is a fibration. ${\mathcal C}_{cf}$ denotes the full subcategory of ${\mathcal C}$ consisting of objects which are simultaneously cofibrant and fibrant.

The following important observation is well-known (\cite[Lemma 4.24]{dwyspa95}, \cite[Proposition 1.2.8]{hov}).

\begin{pro}\label{P:inverting}
Suppose ${\mathcal C}$ is a model category. Then a map of ${\mathcal C}_{cf}$ is a weak equivalence if and only if it is a homotopy equivalence.
\end{pro}


\subsection{Extension Types and Extension Dimension}\label{SS:extensiontypes}
For spaces $X$ and $L$, the notation $L \in {\rm A(N)E}(X)$ means that every map $f \colon A \to L$, defined on a closed subspace $A$ of $X$, admits an extension $\tilde{f} \colon X \to L$ (respectively, $\tilde{f} \colon G \to L$) over $X$ (respectively, over a neighborhood $G$ of $A$ in $X$).

This notation
allows us to define a preorder relation $\leq$ on the class
of CW-complexes. Following \cite{D}, we say that $L \leq K$ if for each space $X$ the condition $L \in {\rm AE}(X)$ implies the condition $K \in {\rm AE}(X)$. The preorder relation $\leq$ naturally generates the equivalence relation: $L\sim K$ iff $L\leq K$ and $K\leq L$. We
denote by $[L]$ the equivalence class of $L$. These equivalence classes of complexes are called {\em extension types}. The above defined relation $\leq$ creates a partial order in the class of extension types. This partial order will still be denoted by $\leq$. Note that under these definitions the class of all extension types has both maximal and minimal elements.
The minimal element is the extension type of the $0$-dimensional sphere $S^{0}$ (i.e. the two-point discrete space) and the maximal element is obviously the extension type of the one-point space $\{ \ast \}$ (or, equivalently, of any contractible complex).

\begin{example}\label{E:order1}
The following observations express some basic properties of the above order.
\begin{itemize}
\item[(i)]
$\min\{ [L], [K]\} = [L\vee K]$, here $L\vee K$ denotes the bouquet of complexes $K$ and $L$.
\item[(ii)]
$[S^{0}] = [L] \Longleftrightarrow L$ is not connected.
\item[(iii)]
$[S^{1}] \leq [L] \Longleftrightarrow L$ is connected.
\item[(iv)]
It follows from the above two observations that there is no complex $L$ such that $[S^{0}] < [L] < [S^{1}]$.
\item[(v)]
Clearly, $S^{n} \in {\rm AE}(X) \Longleftrightarrow \dim X \leq n$.
\item[(vi)]
Similarly, $K(G,n) \in {\rm AE}(X) \Longleftrightarrow \dim_{G}X \leq n$. Here $\dim_{G}X$ stands for the {\em cohomological dimension}\index{cohomological dimension} of $X$ with coefficients in an abelian group $G$ and $K(G,n)$ denotes the corresponding {\em Eilenberg-MacLane complex}\index{Eilenberg-MacLane complex}, i.e. a complex satisfying the following conditions: $\pi_{n}(K(G,n)) = G$ and $\pi_{k}(K(G,n)) = 0$ for each $k \neq n$.
\item[(vii)]
Obviously $[S^{n}] \leq [K(\mbox{\f Z},n)]$, but $[S^{n}] \neq [K(\mbox{\f Z},n)]$. The last part follows from \cite{D1} (for $n=3$) and \cite{dywalsh} (for $n=2$).
\item[(viii)]\label{E:8}
$[S^{n}] < \left[ M({\mathbb Z}_{2},n+1) \vee S^{n+1}\right] < [S^{n+1}]$, where $M({\mathbb Z}_{2},n+1)$ is the corresponding Moore complex, $n \geq 1$.
\item[(ix)]\label{E:9}
The extension type $[{\mathbb R}P^{2}]$ of the projective plane is not comparable with $[S^{n}]$ for any $n \geq 2$ (see \cite{DR}).
\end{itemize}
\end{example}

The Homotopy Extension Theorem implies the following trivial observation.
\begin{pro}\label{P:homotopyequivalent}
If $L$ and $K$ are homotopy equivalent complexes, then $[L] = [K]$.
\end{pro}

Observe that $[S^{n} \vee S^{n+1}] = [S^{n}]$ which shows that homotopy inequivalent complexes might have the same extension type.

The following notion is introduced by A.~Dranishnikov (see \cite{D} and \cite{DD}). The
{\it extension dimension\/} of a~Tychonov space $X$ is less than or
equal to $[L]$ (briefly, $\ed(X)\le[L]$) if $X \in AE(L)$. More precisely
\[\ed(X) = \min\{ [L] \colon L \in {\rm AE}(X)\} .\]


\section{$[L]$-homotopies and $[L]$-complexes}\label{S:lhomotopies}
Throughout this section $L$ stands for a finite CW-complex.
Let $A$ be a subspace of a Polish space $X$ and let $f_0, f_1\colon X \to Y$ be two maps such that $f_{0}(x) = f_{1}(x)$ for each $x \in A$. Then $f_{0}$ and $f_{1}$ are said to be {\it $[L]$--homotopic relative to $A$} (notation: $f_0\stackrel{[L]}\simeq f_1\operatorname{rel} A$) if for any map $h\colon Z\to X\times [0,1]$ where $Z$ is a Polish space of extension dimension $\ed Z\le [L]$, there exists a map $H \colon Z \to Y$ such that
\[ H(z) = \begin{cases}
h(f_{0}(z)),\;\text{if}\; z \in h^{-1}(X\times\{ 0\} \cup A \times [0,1]),\\
h(f_{1}(z)),\;\text{if}\; z \in h^{-1}(X \times \{ 1\} ).
\end{cases}
\]
The following diagram illustrates the situation:

\[
\dgARROWLENGTH=2.0\dgARROWLENGTH
\dgARROWPARTS=8
\begin{diagram}
\node{Z} \arrow{e,l}{h} \arrow{ese,l,6}{H} \node{X\times [0,1]}\\
\node{\widetilde{Z}} \arrow{n,l}{j} \arrow{e,r}{h}
\node{X\times\{ 0,1\} \cup A\times [0,1]} \arrow{n,r,7}{i} \arrow{e,r}{\phi} \node{Y\; .}
\end{diagram}
\]
\medskip

\noindent Here $\widetilde{Z} = h^{-1}\left( X\times\{ 0,1\} \cup A\times [0,1]\right)$ and

\[ \phi (x) =
\begin{cases}
f_{0}(x),\; \text{if}\; x \in X \times \{ 0\} \cup A \times [0,1],\\
f_{1}(x),\; \text{if}\; x \in X \times \{ 1\} .
\end{cases}
\]

If $A =\emptyset$, then we say that $f_{0}$ and $f_{1}$ are $[L]$-homotopic (notation: $f_{0} \stackrel{[L]}\simeq f_{1}$, see \cite[Definition 2.9]{chi}).

It is easy to verify \cite[Proposition 2.7]{chi} that extension dimension of a locally compact polyhedron is identical with its standard dimension. However the concept of $[L]$-homotopy differs from the classical concepts of homotopy or $n$-homotopy even for maps between finite polyhedra -- this can be easily understood analyzing the identity map of the complex indicated in example \ref{E:order1}(viii). 

The class of spaces with respect to which $[L]$--homotopies behave well is identified in the following definition.

\begin{defin}[Definition 2.2, \cite{chi}]\label{D:absoluteextensors}
We say that a space $X$ is an absolute (neighborhood) extensor modulo $L$, or shortly that $X$ is an ${\rm ANE}([L])$-space (notation: $X \in {\rm ANE}([L])$) if  $X \in {\rm ANE}(Y)$ for each space $Y$ with $\ed(Y) \leq [L]$.
\end{defin}

Some of the basic properties of ordinary homotopies have their analogs for $[L]$-homotopies (see \cite{chi} for details). Here is an analog of the standard homotopy extension theorem \cite[Proposition 2.28]{chi}.

\begin{pro}\label{P:globalhomotopyext}
Let $L$ be a finitely dominated complex and $X$ be a Polish ${\rm ANE}([L])$-space. Suppose that $A$ is closed in a space $B$ with $\ed(B) \leq [L]$. If  maps $f , g \colon A \to X$ are $[L]$-homotopic and $f$ admits an extension $F \colon B \to X$, then $g$ also admits an extension $G \colon B \to X$, and it may be assumed that $F
\stackrel{[L]}{\simeq} G$.
\end{pro}

We also recall the following statement \cite[Proposition 2.26]{chi}.

\begin{pro}\label{P:closemaps}
Let $L$ be a finitely dominated complex and $X$ be a Polish ${\rm ANE}([L])$-space $X$. Then there exists an open cover ${\mathcal U} \in \operatorname{cov}(X)$ such that any two $\mathcal U$-close maps of any space into $X$ are $[L]$--homotopic.
\end{pro}

The class of approximately $[L]$-soft maps plays an important role in the $[L]$-homotopy theory. Between locally nice spaces of extension dimension not exceeding $[L]$, such maps provide basic examples of $[L]$-homotopy equivalences.

\begin{defin}[Definition 2.6, \cite{chi}]\label{D:approxsoft}
A map $f \colon X \to Y$ is said to be {\em approximately $[L]$-soft}, if for each space $B$ with $\ed(B) \leq [L]$, for each closed subset $A$ of it, for an open cover ${\mathcal U} \in \operatorname{cov}(Y)$, and for any two maps $g \colon A \to X$ and $h \colon B \to Y$ such that $fg = h/A$, there is a map $k \colon B \to X$ satisfying the conditions $k/A = g$ and the composition $fk$ is ${\mathcal U}$-close to $h$.
\end{defin}

Important examples of approximately $L$-soft maps between polyhedra are those whose non-trivial fibers are topological (or, more generally, homotopical) copies of the complex $L$.

\begin{pro}\label{P:equiv}
Let $p \colon X \to Y$ be an approximately $[L]$-soft map between Polish spaces. Let also
$f_{1}, f_{2} \colon A \to X$ be two maps, defined on a Polish space $A$,  such that
\begin{itemize}
\item[(a)]
$f_{1}(a_{0}) = f_{2}(a_{0})$ for some point $a_{0} \in A$;
\item[(b)]
$p\circ f_{0} \stackrel{[L]}\simeq p\circ f_{1}\operatorname{rel}a_{0}$.
\end{itemize}
Then $f_{0} \stackrel{[L]}\simeq f_{1}\operatorname{rel}a_{0}$.
\end{pro}
\begin{proof}
Consider the map $\phi \colon A \times \{ 0,1\} \cup \{ a_{0}\} \times I\to Y$, defined by letting
\[ \phi (a,t) =
\begin{cases}
p(f_{t}(a)),\; \text{if}\; a \in A\;\text{and}\; t=0,1;\\
p(f_{0}(a_{0})),\; \text{if}\; a = a_{0} \;\text{and}\; t \in I .
\end{cases}
\]

\noindent Let also $h \colon Z \to A \times I$ be an $[L]$-invertible map such that $\ed Z \leq [L]$ and $\widetilde{Z} = h^{-1}\left( A \times\{ 0,1\} \cup \{ a_{0}\} \times I \right)$. By (b), there exists a map $H \colon Z \to Y$ which extends the composition $\phi\circ h |\widetilde{Z} \colon \widetilde{Z} \to Y$.

Next consider the map $\varphi \colon A \times \{ 0,1\} \cup \{ a_{0}\} \times I \to X$

\[ \varphi (a,t) =
\begin{cases}
f_{t}(a),\; \text{if}\; a \in A\;\text{and}\; t=0,1;\\
f_{0}(a_{0}),\; \text{if}\; a = a_{0} \;\text{and}\; t \in [0,1] .
\end{cases}
\]

\noindent which, according to (a), is well defined.

\noindent It is easy to see that the following diagram of unbroken arrows commutes (here $i$ and $j$ denote the corresponding inclusion maps).
\[
\dgARROWPARTS=8
\dgARROWLENGTH=1\dgARROWLENGTH
\begin{diagram}
\node[2]{A\times\{ 0,1\} \cup \{ a_{0}\} \times I} \arrow{e,l}{\varphi} \arrow[3]{s,l}{i} \arrow{k,l}{\phi} \node{X}\arrow[4]{s,l}{p}\\
\node{\widetilde{Z}} \arrow{ne,r}{h} \arrow[3]{s,r}{j}\\
 \\
\node[2]{A\times I}\\
\node{Z} \arrow[2]{nne,r,7,..}{G} \arrow{ne,r}{h} \arrow[2]{e,l}{H} \node[2]{Y}
\end{diagram}
\]

\noindent In particular, $p\circ \varphi \circ h |\widetilde{Z} = \phi\circ h |\widetilde{Z} = H|\widetilde{Z}$. Finally, since $p$ is approximately $[L]$-soft and since $\ed Z \leq \{L]$ it follows that there exists a map $G \colon Z \to X$ (the broken arrow in the above diagram) such that $G|\widetilde{Z} = \varphi\circ h |\widetilde{Z}$ (notice also that above diagram is not commutative if $G$ is included in it, however, it can be made approximately commutative). This proves that $f_{0} \stackrel{[L]}\simeq f_{1}\operatorname{rel}a_{0}$.
\end{proof}


\subsection{$[L]$-complexes}\label{SS:complexes}
It is well known that one of the primary goals of introducing the concept of $CW$-complexes was the possibility of developing a satisfactory homotopy theory. The class of spaces which is specifically designed for the needs of $[L]$-homotopy theory can be defined similarly. First we need the following resolution theorem (statements (i)--(iii) below are contained in [Proposition 2.23]\cite{chi}).

\begin{pro}[Proposition 2.23, \cite{chi}]\label{P:mengerL}
Let $L$ be a finite complex and $X$ be a locally finite polyhedron. Then there exist a locally compact metrizable space $\mu_{X}^{[L]}$ and an $[L]$-invertible and approximately $[L]$-soft proper map $f_{X}^{[L]} \colon \mu^{[L]}_{X} \to X$ satisfying the following conditions:
\begin{itemize}
\item[(i)] $\mu_{X}^{[L]} \in {\rm ANE}([L])$.
\item[(ii)]
$\ed(\mu_{X}^{[L]}) = [L]$.
\item[(iii)] For any map $f \colon B
\to \mu_{X}^{[L]}$, where $B$ is a compact space with $\ed(B) \leq
[L]$, and for any open cover ${\mathcal U} \in
\operatorname{cov}(\mu_{X}^{[L]})$ there is an embedding $g \colon
B \to \mu_{X}^{[L]}$ which is ${\mathcal U}$-close to $f$ and such
that $f_{X}^{[L]}\circ g =f_{X}^{[L]}\circ f$.
\item[(iv)] If, in
addition, $\tau$ is a 
triangulation of $X$, then one can assume that for any subpolyhedron $Y$ of
$X$ (with respect to $\tau$), the
inverse image $\left( f_{X}^{[L]}\right)^{-1}(Y)$ is also a
locally compact ${\rm ANE}([L])$-space and the restriction\\
$f_{X}^{[L]}\left.\right|\left( f_{X}^{[L]}\right)^{-1}(Y)
\colon  \left( f_{X}^{[L]}\right)^{-1}(Y) \to Y$ is also
approximately $[L]$-soft.
\end{itemize}
\end{pro}
\begin{proof}
Proof of this statement, with minor adjustments, follows
Dranishnikov's construction \cite{DR}. Lemmas 2.6 and 2.7 from
\cite{chi} (see also \cite[Lemmas 2.2 and 2.3]{DR})
allow us to construct inductively an inverse sequence
${\mathcal S} = \{ X_{n}, p_{n}^{n+1}\}$, consisting of locally
compact polyhedra $X_{n}$ (with certain triangulations whose
meshes converge to zero) and $[L]$-invertible,
approximately $[L]$-soft, proper, simplicial
bonding maps $p_{n}^{n+1} \colon X_{n+1}\to X_{n}$ so that $X_{1}$
is the given polyhedron
$X$ considered with the given triangulation $\tau$. As in \cite[Lemma 2.3]{DR}, we may assume that ${\mathcal S}$ is $L$-resolvable inverse sequence. We let
$\mu_{X}^{[L]} = \lim{\mathcal S}$ and let the map
$f_{X}^{[L]} \colon \mu^{[L]}_{X} \to X$ be the limit projection
$p_{\infty} \colon \lim{\mathcal S} \to X_{1}$ of the given spectrum.
As in the proof of \cite[Theorem 2.4]{DR}, we have $\ed(\mu_{X}^{[L]}) \leq [L]$ (Property (ii)). Property (i) follows from
\cite[Proposition 2.22]{chi}. Property (iv) is satisfied by construction. If, during the inductive construction, we define $X_{n+1}$ insuring that the bonding map $p_{n}^{n+1} \colon X_{n+1} \to X_{n}$ factors through the projection $X_{n} \times [0,1] \to X_{n}$, then one can easily verify that property (iii) would be also satisfied. 
\end{proof}

Below, for each $n \geq 0$, we consider at most $[L]$-dimensional
compact ${\rm AE}([L])$-spaces $D^{n}_{[L]}$ which admit
($[L]$-invertible) approximately $[L]$-soft maps onto
$n$-dimensional disk $D^{n}$. Actually Proposition
\ref{P:mengerL}(iv) allows us to consider pairs $(S_{[L]}^{n},
D_{[L]}^{n+1})$ as approximately $[L]$-soft preimages of the
standard pairs $(S^{n}, D^{n+1})$ consisting of the
$(n+1)$-dimensional disk and its boundary $\partial D^{n+1} =
S^{n}$, $n \geq 0$. We say that $D_{[L]}^{n}$ is an
$n$-$[L]$--disk and that $S_{[L]}^{n}$ is an $n$-$[L]$--sphere.
Note that generally speaking for any given $n$ there are many
$n$-$[L]$--disks and $n$-$[L]$--spheres. Clearly all
$n$-$[L]$--disks, as ${\rm AE}([L])$-compacta, are $[L]$-homotopy
equivalent to the one point spaces. Proposition \ref{P:equiv}
shows that all $n$-$[L]$--spheres are also $[L]$-homotopy
equivalent.

The following concept has, in fact, been introduced in \cite[Section 2.6.2]{chi}. 

\begin{defin}\label{D:L-complex}
We say that a space $X$ of extension dimension $\ed X\le [L]$ is a (finite)
$[L]$--polyhedron if it admits a proper approximately $[L]$--soft map
$f\colon X\to Y$ onto a locally finite polyhedron $Y$ such that
$f$ possesses property (iv) of Proposition \ref{P:mengerL} for
some (finite) triangulation $\tau$ of $Y$. By $[L]$-complexes we understand spaces $X$ with $\ed X\le [L]$ that are $[L]$--homotopy equivalent to $[L]$-polyhedra.
\end{defin}

The above definition guarantees that any $[L]$--polyhedron can  be
constructed by attaching ``cells" (i.e. $[L]$-disks) along their
``boundaries" (i.e. $[L]$-spheres) by means of inclusion maps.
Therefore a wide variety of properties of $[L]$-polyhedra can be
obtained by straightforward modifications of standard
constructions and proofs (see, for instance, \cite{fritch},
\cite{switzer}).

Next, suppose that $X$ is an $[L]$--polyhedron and $f\colon X\to Y$
is a corresponding approximately $[L]$--soft mapping onto a locally
finite polyhedron $Y$. Let $Y^{(n)}$ denote the $n$--skeleton of
$Y$. Then for any $y\in Y^{(0)}$ we have $f^{-1}(y)\in AE([L])$.
Denote $X^{(n)}= f^{-1} (Y^{(n)})$. We call $X^{(n)}$ {\it a
$n$-$[L]$--skeleton} of $X$. Clearly, $X^{(n)}$ is closed in $X$
for any $n$. Since $Y$ is locally finite for any point $x\in X$
there exists an open (in $X$) neighbourhood $Ox$ such that
$Ox\subset X^{(n)}$ for some $n$. This implies that a subset $F$
of $X$ is closed if and only if $F\cap X^{(n)}$ is closed for any
$n$.

Every Hilbert cube manifold, according to the corresponding triangulation theorem, is topologically homeomorphic to the product $K \times {\mathbb Q}$ of a locally compact polyhedron $K$ and the Hilbert cube ${\mathbb Q}$. Since the projection $K \times {\mathbb Q} \to K$ is a proper (approximately) soft map, it follows that every Hilbert cube manifold is a $[\{ \ast\}]$-polyhedron. It is important to note however that every $[\{ \ast\}]$-polyhedron is homotopically equivalent to a standard polyhedron. 

$\mu^{n}$-manifolds (i.e. spaces locally homeomorphic to the $n$-dimensional universal Menger compactum $\mu^{n}$) can also be ``triangulated" in certain non-standard sense (see a comprehensive discussion of related matters in \cite{DR1}, \cite[Chapter 4]{chibook}). In particular, they admit proper approximately $[S^{n}]$-soft maps onto standard polyhedra and consequently are $[S^{n}]$-polyhedra in our sense. They are $[S^{n}]$-homotopy equivalent to standard $n$-dimensional polyhedra.


\subsection{$[L]$-homotopy groups}\label{SS:groups}
Let $S^n$ denote a unit $n$--sphere. Fix a point $s\in S^n$. For
each $n\ge 0$  we consider an
{\it $n$-$[L]$--sphere}, i.e. a compactum of extension dimension
at most $[L]$ which admits an approximately $[L]$--soft mapping onto $S^n$.

Let $(X,x_0)$ be a pointed space.
Let $\sfr{n}$ be a $n$-$[L]$--sphere, $n\ge 1$, and
$f\colon\sfr{n}\to S^n$ be an approximately $[L]$--soft
mapping. Fix a point $\pt\in f^{-1} (s)$ and consider the
set $\grup{n} (X,x_0 )= \left[ (\sfr{n} ,\pt ),(X,x_0 )\right]_{[L]}$
of relative $[L]$--homotopy classes of maps of pointed spaces.

Next, we show that this set does not depend on the choice of $[L]$--sphere.
Consider another $n$-$[L]$--sphere $Q_{[L]}^{n}$ and let
$g \colon Q_{[L]}^{n}\to S^n$ be the corresponding approximately
$[L]$--soft mapping. Choose a point $q_{[L]}\in g^{-1} (s)$ and let
$\widetilde{\pi}_{n}^{[L]}(X,x_0 )=\left[ (Q_{[L]}^{n} ,q_{[L]}),(X,x_0 )\right]_{[L]}$.

Since $S^{n}$ is an ${\rm ANE}$-compactum there exists an
open cover ${\mathcal U} \in \operatorname{cov}(S^{n})$ such
that any two ${\mathcal U}$-close maps (defined on any compactum)
are homotopic as maps into $S^{n}$
(this is Proposition \ref{P:closemaps} with $L = \{ \ast\}$).

Now consider the following commutative diagram (consisting of unbroken arrows)

\[
\dgARROWPARTS=8
\dgARROWLENGTH=1\dgARROWLENGTH
\begin{diagram}
\node{\{ s_{[L]}\} } \arrow{e,l}{\alpha} \arrow{s,l}{i} \node{ Q^{n}_{[L]}} \arrow{s,l}{g}\\
\node{S^{n}_{[L]}} \arrow{e,r}{f} \arrow{ne,r,..}{p } \node{S^{n}}\\
\end{diagram}
\]

\noindent where $\alpha (s_{[L]}) =  q_{[L]}$.

Since $g$ is approximately $[L]$--soft, there exists a
mapping $p\colon \sfr{n}\to Q_{[L]}^{n}$ (the broken arrow
in the above diagram)
such that $p(\pt )= q_{[L]}$ (i.e. $p\circ i = \alpha$) and
$g\circ p$ is ${\mathcal U}$-close to $f$.
Similarly, there exists a mapping $q \colon Q_{[L]}^{n}\to\sfr{n}$
such that $q(q_{[L]})=\pt$ and $f\circ q$ is ${\mathcal U}$-close to $g$.
Choice of the cover ${\mathcal U}$ guarantees that
$g\circ p \simeq f$ and $f\circ q \simeq g$. We may assume
without loss of generality that these are homotopies relative
to the given points $s_{[L]}$ and $q_{[L]}$.

Next note that

\[ g \circ p \circ q \simeq f \circ q \operatorname{rel}q_{[L]} \simeq g \operatorname{rel}q_{[L]} \]

\noindent and
\[ f \circ q \circ p \simeq g \circ p \operatorname{rel}s_{[L]} \simeq f \operatorname{rel}s_{[L]}.\]

According to Proposition \ref{P:equiv}, $\displaystyle q \circ p \stackrel{[L]}\simeq \operatorname{id}_{Q^{n}_{[L]}}\operatorname{rel}q_{[L]}$
and $\displaystyle p \circ q \stackrel{[L]}\simeq \operatorname{id}_{S^{n}_{[L]}}\operatorname{rel}s_{[L]}$. We shall
refer to mappings constructed as described above as {\it
$[L]$--homotopy equivalences of canonical type}. Observe, that any two $[L]$--homotopy equivalences of canonical type are $[L]$--homotopic.

Define a mapping $\phi\colon\grup{n} (X,x_0 )\to \widetilde{\pi}_{n}^{[L]}(X,x_0 )$ as follows. Consider an element $\alpha\in\grup{n} (X,x_0)$ and let $a\colon (\sfr{n},\pt )\to (X,x_0 )$ be its representative. We let $\phi (\alpha )=\beta$, where $\beta$ is a
 relative $[L]$--homotopy class of the composition $a\circ q\colon (Q_{[L]}^{n}, q_{[L]})\to (X,x_0 )$. Similarly, we define a mapping $\psi\colon \widetilde{\pi}_{n}^{[L]}(X,x_0 )\to\grup{n} (X,x_0 )$ using the mapping $p$.
It is easy to check that $\phi$ and $\psi$ are well-defined. Clearly they are inverses to each other. Thus, there exists a bijection between the sets $\grup{n} (X,x_0 )$ and $\widetilde{\pi}_{n}^{[L]}(X,x_0 )$.
Therefore, the set $\grup{n} (X,x_0)$ does not depend on the choice of $n$-$[L]$--sphere.
It is clear that bijections $\phi$ and $\psi$ do not depend on the choice of $[L]$--homotopy equivalences of canonical type.

Let $f\colon (X,x_0)\to (Y,y_0)$ be a mapping and $[f]_{[L]}\in [(X,x_0),(Y,y_0)]_{[L]}$ be its relative $[L]$--homotopy class. Then a natural map
$\grup{n} ([f]_{[L]})\colon\grup{n} (X,x_0)\to\grup{n} (Y,y_0)$ can be defined in a standard way.

Our next goal is to introduce a group structure on $\grup{n} (X,x_0)$ such that natural maps
$\grup{n} ([f]_{[L]})$ are well-defined homomorphisms. Moreover, this
structure will be defined so that
bijections $\phi$ and $\psi$, generated by $[L]$--homotopy
equivalences of canonical type, are group isomorphisms.

Let $\alpha$ and $\beta$ be two elements of
the set $\grup{n} (X,x_0 )$ and $a,b\colon (\sfr{n},
\pt )\to (X,x_0)$ be their representatives. Let $S^n_+$ and $S^n_-$ denote the upper and lower hemispheres, respectively, and $E$ denote an equator of $S^n$ containing the point $s$. Let $f\colon\sfr{n}\to S^n$ be an approximately $[L]$--soft mapping.
Let $h\colon S^n\to S^{n}\vee S^{n}$ be the homotopy comultiplication defining the standard $H$-cogroup structure (see, for instance, \cite[Definition 2.16]{switzer}) of the sphere $S^{n}$. Let also $f_{-} \colon S^{n}_{[L]} \to S^{n}$ and $f_{+} \colon S^{n}_{[L]} \to S^{n}$ be two copies of the map $f$.

Since $f_{+}$ is approximately $[L]$-soft there exists a mapping
$\widetilde{a}\colon f^{-1}(S^{n}_{+}) \to \sfr{n}$ such that $\widetilde{a}(f^{-1}(E))= s_{[L]}$ and the composition $f_{+}\circ \widetilde{a}$ is
${\mathcal U}$--close to the composition $h\circ f|_{f^{-1}(S^{n}_{+})}$.
Similarly, there exists a mapping $\widetilde{b}\colon f^{-1}(S^{n}_{-})\to \sfr{n}$ such that $\widetilde{b}(f^{-1}(E))= s_{[L]}$ and the composition $f_{-}\circ \widetilde{b}$ is
${\mathcal U}$--close to the composition $h\circ f|_{f^{-1}(S^{n}_{-})}$.
Here is the corresponding diagram (commutative up to $[L]$--homotopy):

\[
\dgARROWPARTS=8
\dgARROWLENGTH=1\dgARROWLENGTH
\begin{diagram}
\node[3]{(X,x_{0})}\\
\node{S_{[L]}^{n} = f^{-1}(S_{-}^{n})\cup f^{-1}(S_{+}^{n})} \arrow{e,l}{\widetilde{a}\cup\widetilde{b}} \arrow{s,l}{f} \arrow{ene,l,..}{a\ast b} \node{S^{n}_{[L]} \vee S^{n}_{[L]}} \arrow{ne,r}{a\vee b} \arrow{s,r}{f_{-}\vee f_{+}}\\
\node{S^{n}} \arrow{e,t}{h}  \node{S^{n}\vee S^{n}}\\
\end{diagram}
\]

\noindent It is routine to check that the relative $[L]$--homotopy class $[a*b]_{[L]}$ of the composition $a\ast b = (\widetilde{a}\cup\widetilde{b})\circ (a\vee b) \colon (\sfr{n},E')\to (X,x_0)$, which formally is defined by letting

\[
 (a\ast b)(x)=\cases
a(\widetilde{a}(x)),
         &\text{if} \; x\in f^{-1}(S_{+}^{n})\\
        b(\widetilde{b}(x)),
         &\text{if} \; x\in f^{-1}(S_{-}^{n}) ,\endcases
\]

\noindent does not depend on the choice of representatives $a$ and $b$ (and mappings $\widetilde{a}$ and $\widetilde{b}$).
Now we can define the product of $\alpha$ and $\beta$ by letting $\alpha\ast\beta = [a*b]_{[L]}$.

The unit element in $\grup{n} (X,x_0 )$ is given by $\epsilon = [e]_{[L]}$ where $e$ is a constant mapping which sends $\sfr{n}$ to the point $x_0$.

Finally, given an element $\alpha\in\grup{n} (X,x_0 )$ we define its inverse $\alpha ^{-1}\in\grup{n} (X,x_0 )$ with respect to the operation $\ast$ as follows.
Let $a\colon (\sfr{n}, \pt )\to (X,x_0)$ be a representative of $\alpha$.
Let $g\colon S^n\to S^n$ be the mapping such that
$g(x_1,\dots ,x_{n+1})= (x_1,\dots ,-x_{n+1} )$, which fixes the equator $E$.
Since $f$ is approximately $[L]$--soft there exists a mapping $\widetilde{g}\colon\sfr{n}\to\sfr{n}$ such that $\widetilde{g}(s_{[L]} )= s_{[L]}$ and composition $f\circ \widetilde{g}$ is ${\mathcal U}$--close to the composition $g\circ f$. Here is the diagram (as before commutative up to $[L]$-homotopy):

\[
\dgARROWPARTS=8
\dgARROWLENGTH=1\dgARROWLENGTH
\begin{diagram}
\node[3]{(X,x_{0})}\\
\node{S_{[L]}^{n}} \arrow{e,l}{\widetilde{g}} \arrow{s,l}{f} \arrow{ene,l,..}{\widetilde{a}} \node{S^{n}_{[L]}} \arrow{ne,r}{a} \arrow{s,r}{f}\\
\node{S^{n}} \arrow{e,t}{g}  \node{S^{n}}\\
\end{diagram}
\]

\noindent The $[L]$-homotopy class of the composition $\widetilde{a} = a\circ \widetilde{g} \colon (S_{[L]}^{n},s_{[L]}) \to (X,x_{0})$ does not depend on the choice of representative $a$ and of the mapping $\widetilde{g}$. This allows us to define $\alpha ^{-1} = [a\circ \widetilde{g}]_{[L]}$. It only remains to note that $\alpha ^{-1}\alpha =
\alpha\alpha ^{-1} =\epsilon\in\grup{n} (X,x_0)$.


\section{Model category structure on
${\mathbf T}{\mathbf O}{\mathbf P}$ generated by a finite complex}\label{S:model}

In this section we prove our main results which states that the category ${\mathbf T}{\mathbf O}{\mathbf P}$ admits a model category structure whose weak equivalences are weak $[L]$-homotopy equivalences. We begin by introducing the needed classes of morphisms.

\begin{defin}\label{D:weak}
A map $f \colon X \to Y$ of spaces is called a {\em weak $[L]$-homotopy
equivalence} if for each basepoint $x \in X$ the map
$f_{\ast} \colon \pi_{n}^{[L]}(X,x) \to \pi_{n}^{[L]}(Y,f(x))$ is a
bijection of pointed sets for $n=0$ and an isomorphism of groups
for $n \geq 1$.
\end{defin}

In order to define cofibrations and fibrations first recall that
a map $i \colon  A \to B$ has the {\em left lifting property} (LLP) with respect
to a map $p \colon X \to Y$ and $p$ has the {\em right lifting
property} (RLP) with respect to $i$ if, for every commutative square diagram of
unbroken arrows

\[
\begin{diagram}
\node{A} \arrow{e,t}{f} \arrow{s,l}{i}  \node{X} \arrow{s,r}{p}\\
\node{B} \arrow{e,t}{g} \arrow{ne,t,..}{h} \node{Y}
\end{diagram}
\]

\noindent there exists a lift (the broken arrow) $h \colon B \to X$ such that
$h\circ i = f$ and $p\circ h = g$.

\begin{defin}\label{D:Lmodel}
Let $f \colon X \to Y$ be a map in
${\mathbf T}{\mathbf O}{\mathbf P}$. We say that $f$ is
\begin{itemize}
\item[(i)$_{L}$]
a weak equivalence if it is a weak $[L]$-homotopy equivalence,
\item[(ii)$_{L}$]
a fibration if it has the RLP with respect to inclusions of
finite $[L]$-polyhedra $A\hookrightarrow B$ such that both $A$ and $B$ are $AE([L])$--spaces.
\item[(iii)$_{L}$]
a cofibration if it has the LLP with respect to acyclic fibrations (see Definition \ref{D:model}).
\end{itemize}
\end{defin}

Below we shall use the following fact, proof of which is trivial.

\begin{pro}\label{P:triv}
Let $f\colon S^n _{[L]}\to X$ be a mapping of $n$-$[L]$--sphere to a topological space $X$.
If $f$ is $[L]$--homotopic to a constant map then for any embedding of $S^n _{[L]}$ to a space $Z$ of extension dimension $\le [L]$ the mapping $f$ can be extended to a mapping of $Z$ to $X$.
Conversely, suppose that $S^n _{[L]}$ is a subspace of a space $Z$ of extension dimension $\le [L]$ such that $Z\in AE([L])$ and $f$ admits an extension over $Z$. Then $f$ is $[L]$-homotopic to a constant map.
\end{pro}

In order to verify axioms of model category we need to obtain a
characterization  of acyclic $[L]$--fibrations in terms of the RLP.

\begin{lem}\label{L:acyc}
An $[L]$-fibration $p\colon X\to Y$ is acyclic iff it has the RLP with
respect to any pair $\sfr{n}\subset\disk{n+1}$, where $\disk{n+1}$
is an $(n+1)$-$[L]$--disk and $\sfr{n}$ is a corresponding
$n$-$[L]$--sphere, $n=0,1,\dots$.
\end{lem}
\begin{proof}
First let us prove the following claim.

{\bf Claim}. {\em Let $p\colon X\to Y$ be an acyclic $[L]$--fibration. Then there
exist pairs of $[L]$-polyhedra $(\sfr{n} ,\disk{n+1} )$, where
$\disk{n+1}$ is an $(n+1)$-$[L]$--disk and $\sfr{n}$ is the
corresponding $n$-$[L]$--sphere, $n=0,1,\dots$, such that $p$ has
the RLP with respect to inclusions $\sfr{n}\subset\disk{n+1}$.}

Let $D^{n+1} _r$ and $S^n _r =\partial D^{n+1} _r$ denote centered
at the origin $O$ (of the $(n+1)$-dimensional Euclidean space
$\mathbb{R}^{n+1}$) $(n+1)$--disk and $n$--sphere of radius $r$,
respectively. Let $f\colon D\to D^{n+1}_1 =D^{n+1}$ be an
approximately $[L]$--soft mapping of a compactum $D$ having
extension dimension $\le [L]$ onto $D^{n+1}$. Put
$\sfr{n}=f^{-1}(S^n _1 )$, $\sfr{n} (0)=f^{-1} (S^n_{1/2})$,
$\cyl{n+1} =f^{-1} (D^{n+1}\backslash\Int (D^{n+1} _{1/2}))$ and
$D_0 =f^{-1} (D^{n+1} _{1/2})$. Let $\tau$ denote a triangulation
on $D^{n+1}$ such that $S^n_{1/2}$, $S^n_1$, $D^{n+1}_{1/2}$ and
$D^{n+1}\backslash\Int (D^{n+1}_{1/2})$ are subpolyhedra of
$\tau$. According to Proposition \ref{P:mengerL}, one can
construct mapping $f$ so that restrictions of $f$ on preimages of
subpolyhedra of $D^{n+1}$ with respect to $\tau$ are also
approximately $[L]$-soft.

This implies, in particular, that $D_0\in AE([L])$, $\sfr{n}\in ANE([L])$ and that $\cyl{n+1}$ is an $(n+1)$-$[L]$--cylinder, i.e. ANE([L])--compactum of extension dimension $\le [L]$ admitting an approximately $[L]$-soft mapping $f\colon\cyl{n+1}\to S^n\times I$, where  $I=[0,1]$ denotes unit interval.

Let $\disk{n+1} = D/D_0$ be a quotient space and $\phi\colon D\to
D/D_0$ be a quotient mapping. Since $D\in AE([L])$ and $D_0\in
AE([L])$ it follows that $\disk{n+1}\in AE([L])$. Put $\pt =\phi
(D_0)$. Define a mapping $g\colon \disk{n+1}\to D^{n+1}$ as
follows. We let $g(\pt )=O$ and for each point $x\in\disk{n+1}$,
distinct from $\pt$, we let $g(x)=h\circ f (x)$, where $h\colon
D^{n+1}\to D^{n+1}$ is a mapping which collapses $D^{n+1} _{1/2}$
to the point $O$. Let also $\tau '$ denote a triangulation on
$D^{n+1}$ obtained from $\tau$ by means of collapsing $D^{n+1}$ to
point (and enlarging of the resulting cellular structure). It is
easy to check that the mapping $g$ is approximately $[L]$-soft and
restriction of $g$ on a preimage of any subpolyhedron of $D^{n+1}$
with respect to $\tau '$ is also approximately $[L]$-soft.
Therefore the constructed space $\disk{n+1}$ is an
$[L]$-polyhedron. Further in the proof of Claim we shall refer to the compacta
constructed above simply as $[L]$-disk, corresponding $[L]$-sphere
and $[L]$-cylinder, and denote them by $\disk{n+1}$, $\sfr{n}$ and
$\cyl{n+1}$, respectively. Notice also that $\disk{n+1}$ is
homeomorphic to a quotient space $\cyl{n+1} /\sfr{n}(0)$.

Let us show that the above constructed pair $(\sfr{n}
,\disk{n+1})$ satisfies condition of the Claim. Consider two
mappings $f\colon\disk{n+1}\to Y$ and $G\colon\sfr{n}\to X$ such
that $p\circ G=f|_{\sfr{n}}$. Since $p$ is acyclic and
$f|_{\sfr{n}}$ represents  trivial element of the group $\grup{n}
(Y)$, there exists a mapping $\overline{G}\colon\disk{n+1}\to X$
extending $G$.

Consider a compactum $\widehat{\sfr{n+1}}= (\disk{n+1})_{+}\cup (\disk{n+1})_{-}$ obtained by gluing of two copies of $\disk{n+1}$ along $\sfr{n}$.
Clearly, $\widehat{\sfr{n+1}}$ admits an approximately $[L]$--soft mapping onto $S^{n+1}$ and hence represents an $(n+1)$-$[L]$--sphere. Fix the point $\point\in (\disk{n+1})_{-}$ (see the construction of $\disk{n+1}$ above).
Define a mapping $h\colon (\widehat{\sfr{n+1}} ,\point )\to (Y,h(\point ))$ as follows:
\[
 h(x)=\cases f(x),
         &\text{if \;$x\in (\disk{n+1})_{+}$}\\
        p\circ\overline{G} (x),
         &\text{if \;$x\in (\disk{n+1})_{-}$}\endcases
\]
Denote $y_0=h(\point )$.
Then $h$ represents an element $\alpha$ of the group $\grup{n+1} (Y,y_0)$.
Since $p$ is acyclic, there exists a mapping $H\colon (\widehat{\sfr{n+1}}
,\point )\to (X,x_0)$ where $x_0\in p^{-1} (y_0)$ such that composition
$p\circ H\colon (\widehat{\sfr{n+1}} ,\point )\to (Y,y_0)$ represents an
element of the group $\grup{n+1} (Y,y_0 )$ which is $[L]$--homotopically
inverse to $\alpha$.

Let $\cyl{n+1}$ be an $n$-$[L]$--cylinder corresponding to $\disk{n+
1}$. Consider the space $\widehat{\disk{n+1}} =\disk{n+1}\cup\cyl{n+1}$
obtained by gluing of $\disk{n+1}$ and $\cyl{n+1}$ along
$\sfr{n}$ and the space
$\widetilde{\sfr{n+1}} =(\widehat{\disk{n+1}} )_{+}\cup (\widehat{\disk{n+1}})_{-}$ obtained by
gluing of two copies of $\widehat{\disk{n+1}}$ along $\sfr{n} (0)$.
Denote by $\phi$ the corresponding quotient mapping. It is easy to check that compactum
$\widetilde{\sfr{n+1}}$ represents a $(n+1)$-$[L]$--sphere.
Let
$$\widetilde{\disk{n+1}}=(\cyl{n+1} )_{+}\cup _{\phi} (\widehat{\disk{n+1}}
)_{-}\subset\widetilde{\sfr{n+1}}$$
\noindent Clearly, $\widetilde{\disk{n+1}}$ is an $(n+1)$-$[L]$--disk and an $[L]$-polyhedron.
In particular, $\widetilde{\disk{n+1}}\in AE([L])$
Define a mapping
$$k\colon (\widetilde{\sfr{n+1}} , \sfr{n} (0))\to (Y,y_0)$$
\noindent as follows:
\[
 k(x)=\cases h(x),
         &\text{if \;$x\in (\widehat{\disk{n+1}})_{+}\backslash\sfr{n} (0)$}\\
         p\circ H (x),
         &\text{if \;$x\in (\widehat{\disk{n+1}})_{-}\backslash\sfr{n} (0)$}\\
         y_0,
         &\text{if \;$x\in\sfr{n} (0)$}\endcases
\]

\noindent Define a mapping $K\colon (\widetilde{\disk{n+1}} )\to X$ such that
$k|_{\widetilde{\disk{n+1}}} = p\circ K$ by letting

\[
 K(x)=\cases \overline{G} (x),
         &\text{if \;$x\in (\cyl{n+1})_{+}\backslash\sfr{n} (0)$}\\
         H (x),
         &\text{if \;$x\in (\widehat{\disk{n+1}})_{-}\backslash\sfr{n} (0)$}\\
         x_0,
         &\text{if \;$x\in\sfr{n} (0)$}\endcases
\]

\noindent The mapping $k$ represents the product of $\alpha$ and $\alpha ^{-1}$ in the group
$\grup{n+1} (Y,y_0)$ and hence is $[L]$--homotopically trivial.

Let $\mu ^{[L]}$ be a strongly $[L]$--universal ${\rm AE}([L])$--compactum\footnote{Note that the Hilbert cube ${\mathbb Q}$ is an example of strongly $[\{\ast\} ]$-universal ${\rm AE}([\{ \ast\} ])$-compactum of extension dimension $[\{ \ast\} ]$ and the universal $n$-dimensional Menger compactum $\mu^{n}$ serves as an example of strongly $[S^{n}]$-universal ${\rm AE}([S^{n}])$-compactum of extension dimension $[S^{n}]$ (see \cite[Chapters 2 \& 4]{chibook} for a comprehensive discussion of various aspects of strong universality property within the general theory of absolute extensors).} of extension dimension $[L]$ provided by Proposition \ref{P:mengerL}. According to (iv) of that proposition we may assume that $\mu ^{[L]}$ is an $[L]$-polyhedron.
Let $i\colon\widetilde{\sfr{n+1}}\hookrightarrow\mu ^{[L]}$ be an embedding. 
By Proposition \ref{P:triv}, there exists an extension $\overline{k}\colon\mu ^{[L]}\to Y$ of $k$.
Applying the RLP of $p$ to the pair $\widetilde{\disk{n+1}}\subset\mu ^{[L]}$ we obtain a lifting
$\overline{K}\colon Z\to X$
of $\overline{k}$ which is an extension of $K$.
The mapping $F =\overline{K} |_{(\disk{n+1})_{+}}\colon
(\disk{n+1})_{+}\to X$ provides a desired extension of $G$.
This concludes the proof of the claim.

Let now $\sfr{n}\subset\disk{n+1}$ be a pair of $[L]$-polyhedra
($(n+1)$-$[L]$--disk and corresponding $n$-$[L]$--sphere) provided
by the Claim. Let also $\tau$ denote the
correspondent triangulation of $D^{n+1}$ and
$f\colon\disk{n+1}\to D^{n+1}$ be the correspondent approximately
$[L]$-soft mapping. Consider $[L]$-dimensional $AE([L])$--space
$\ddisk{n+1}$ which satisfies conditions (i)-(iv) of Proposition
\ref{P:mengerL} and the corresponding $[L]$-invertible and
approximately $[L]$-soft mapping $\fL\colon\ddisk{n+1}\to
D^{n+1}$. We assume that condition (iv) of Proposition
\ref{P:mengerL} is satisfied with respect to the same
triangulation $\tau$ of $D^{n+1}$. Put $\ssfr{n} =(\fL )^{-1} (S^n )$.
Then $\ssfr{n}\in ANE([L])$ and the restriction $\fL |_{\ssfr{n}}$
is also approximately $[L]$-soft.

The further proof consists of two steps. First we apply the RLP of
$p$ with respect to the pair $\sfr{n}\subset\disk{n+1}$ to verify
that $p$ has RLP with respect to the pair $(\ssfr{n} ,\ddisk{n+1}
)$. Then we use this fact to show that $p$ has RLP with respect to
arbitrary pair of $n$-$[L]$--sphere and $(n+1)$-$[L]$--disk.

Since the mapping $\fL$ is $[L]$-invertible there exists a mapping
$\tilde{f}$ such that $\fL\circ\tilde{f} =f$. Next, Proposition
\ref{P:mengerL}(iii) provides us with an embedding
$g\colon\disk{n+1}\to\ddisk{n+1}$ such that
$$\fL\circ g=\fL\circ\tilde{f} =f.$$
It is easy to check that
$$g(\sfr{n} )\subset\ssfr{n}\mbox{ and
}g(\disk{n+1})\backslash\sfr{n}\subset\ddisk{n+1}\backslash\ssfr{n}.$$
In what follows we identify compacta $\disk{n+1}$ and $\sfr{n}$
with their images $g(\disk{n+1})$ and $g(\sfr{n})$. Consider the
space $\tdisk{n+1}=\ssfr{n}\cup\disk{n+1}$. It is not hard to see
that  $\tdisk{n+1}\in AE([L])$. Moreover,
$\tdisk{n+1}$ is $(n+1)$-$[L]$--disk and $[L]$-polyhedron.
Let now
$F\colon\ddisk{n+1}\to Y$ and $G\colon\ssfr{n}\to X$ be mappings
such that $F|_{\ssfr{n}} =p\circ G$. We have to find a lifting
$\overline{G}$ of $F$ extending $G$. Since $p$ has RLP with respect
to the pair $\sfr{n}\subset\disk{n+1}$ there exists an extension
$\widetilde{G}\colon\tdisk{n+1}\to X$ of $G$ such that
$F|_{\tdisk{n+1}} =p\circ\widetilde{G}$. Now the right lifting
property of $[L]$-fibration applied to the pair
$\tdisk{n+1}\subset\ddisk{n+1}$ provides us with the desired lifting
$\overline{G}$.

Let now $(\sfr{n} ,\disk{n+1})$ denote an arbitrary pair of
$(n+1)$-$[L]$--disk and $n$-$[L]$--sphere. The corresponding
approximately $[L]$ soft mapping of $\disk{n+1}$ onto $D^{n+1}$ is
denoted still by $f$. As above, we can find an embedding
$\disk{n+1}\hookrightarrow\ddisk{n+1}$ such that $\sfr{n}
\hookrightarrow\ssfr{n}$ and
$\disk{n+1}\backslash\sfr{n}\hookrightarrow\ddisk{n+1}\backslash\ssfr{n}.$
Since both mappings $\fL |_{\ssfr{n}}$ and $f|_{\sfr{n}}$ are
approximately $[L]$--soft there exists a retraction
$r\colon\ssfr{n}\to\sfr{n}$ (which is $[L]$-homotopy equivalence).
Since $\disk{n+1}\in AE([L])$ this retraction can be extended to a
retraction $R\colon\ddisk{n+1}\to\disk{n+1}.$ Consider two
mappings $F\colon\disk{n+1}\to Y$ and $G\colon\sfr{n}\to X$ such
that $F|_{\sfr{n}} =p\circ G$. Since, as shown above, the mapping
$p$ possesses RLP with respect to the pair
$\ssfr{n}\subset\ddisk{n+1}$ the mapping $F\circ
R\colon\ddisk{n+1}\to Y$ has a lifting
$\overline{G}\colon\ddisk{n+1}\to X$ which extends the mapping
$G\circ r\colon\ssfr{n}\to X$. Clearly, the restriction
$\overline{G} |_{\disk{n}}\to X$ is a lifting of $F$ extending
$G$, as required. Proof of the necessity is completed.

Now we show that the condition of the Lemma is also sufficient.
Let $p\colon X\to Y$ be a fibration possessing the RLP with
respect to any pair $\sfr{n}\subset\disk{n+1}$ where $\disk{n+1}$
is an $(n+1)$-$[L]$--disk and $\sfr{n}$ is a corresponding
$n$-$[L]$--sphere, $n=0,1,\dots$. We need to show that $p$ is
acyclic. Fix $n=0,1,\dots$ and a point $x_0\in X$. Let $y_0 =p(x_0
)$. Consider an arbitrary element $\alpha$ of the group $\grup{n}
(Y,y_0 )$. Let $\widetilde{\sfr{n}} = \disk{n}/\sfr{n-1}$ be a
quotient space, where $\disk{n}$ is an $n$-$[L]$--disk. Let $\phi$
be the correspondent quotient mapping and denote $\point =\phi
(\sfr{n-1} )$. It is easy to check that compactum
$\widetilde{\sfr{n}}$ admits an approximately $[L]$--soft mapping
onto $S^n$. Therefore element $\alpha$ can be represented by means
of a mapping $\widetilde{f}\colon \widetilde{\sfr{n}}\to Y$ such
that $\widetilde{f}(\point )=y_0$. The mapping $\widetilde{f}$
allows us to define mapping $f\colon\disk{n+1}\to Y$ as follows.
\[
 f(x)=\cases \widetilde{f}(x),
         &\text{if \;$x\in\disk{n}\backslash\sfr{n-1}$}\\
         x_0,
         &\text{if \;$x\in\sfr{n-1}$}\endcases
\]

Let $G\colon\sfr{n-1}\to x_0$ be a constant mapping. Applying RLP
of $p$ to the pair $\sfr{n-1}\subset\disk{n}$ we can find a
lifting $F\colon\disk{n}\to X$ of $f$ extending $G$. Now we define
mapping $\widetilde{F}\colon \widetilde{\sfr{n}}\to X$ letting
\[
 \widetilde{F}(x)=\cases F(x),
         &\text{if \;$x\in \widetilde{\sfr{n}}\backslash\point$}\\
         x_0,
         &\text{if \;$x=\point$}\endcases
\]

It is easy to see that $\widetilde{f}=p\circ\widetilde{F}$. This
shows that $p$ induces an epimorphism of $n^{th}$-$[L]$--homotopy
groups.

Consider now a mapping $F\colon\sfr{n}\to X$ where $\sfr{n}$
coincides with the subset $\sfr{n}$ of $(n+1)$-$[L]$--disk
$\disk{n+1}$, such that a composition $p\circ F\colon\sfr{n}\to Y$
is $[L]$--homotopically trivial. Then there exists an extension
$g\colon\disk{n+1}\to Y$ of $p\circ F$. Applying RLP of $p$ to the
pair $\sfr{n}\subset\disk{n+1}$ we can find a lifting
$G\colon\disk{n+1}\to X$ of $g$ extending $F$. By Proposition
\ref{P:triv} the mapping $F$ is also $[L]$--homotopically trivial.
This shows that $p$ induces a monomorphism of the 
$n^{th}$-$[L]$--homotopy groups. Therefore the fibration $p$ is acyclic.
\end{proof}

The above Lemma allows us to prove the following statement.

\begin{cor}\label{C:cofibrant}
Every $[L]$-polyhedron is cofibrant and fibrant object in the category ${\mathbf T}{\mathbf O}{\mathbf P}$ with weak equivalences,
fibrations and cofibrations, defined as in Definition \ref{D:Lmodel}.
\end{cor}
\begin{proof}
Definition of $[L]$--fibration implies that every topological space is a fibrant object. To prove that
every $[L]$-polyhedron $X$ is cofibrant we need to show that for any
acyclic fibration $p\colon A\to B$ and for any mapping $f\colon
X\to B$ there exists a lifting $F\colon X\to A$. Since any
$[L]$--polyhedron can be obtained by attaching $n$-$[L]$-cells, such a
lifting $F$ can be constructed by induction. If $g\colon X\to Y$
is an approximately $[L]$-soft mapping corresponding to $X$ (see
Definition \ref{D:L-complex}), then for each $y\in Y^{(0)}$ we
have $g^{-1}(y)\in AE([L])$. Therefore we can begin the inductive
construction applying RLP of $p$ to the pairs $x\in g^{-1}(y)$ for
each $y\in Y^{(0)}$ (the point $x$ is arbitrary). We perform
inductive steps applying Lemma \ref{L:acyc}. Finally, the
resulting map $F$ is continuous by virtue of the fact that $C$ is
closed subspace of $X$ if and only if $C\cap X^{(n)}$ is closed
for each $n$ (see remarks following Definition \ref{D:L-complex}).
\end{proof}

Next we need the following proposition.

\begin{lem}\label{L:factor}
Every map $p \colon X \to Y$ in ${\mathbf T}{\mathbf O}{\mathbf P}$
can be factored in either of two ways:
\begin{itemize}
\item[(a)]
$f=p_{\infty}i_{\infty}$, where
$i_{\infty} \colon X \to X^{\prime}$ is a cofibration
and $p_{\infty} \colon X^{\prime} \to Y$ is an acyclic $[L]$-fibration.
\item[(b)]
$f=q_{\infty}j_{\infty}$, where
$j_{\infty} \colon X \to X^{\prime}$ is a weak $[L]$-homotopy
equivalence which has the LLP with respect to $[L]$-fibrations,
and $q_{\infty} \colon X^{\prime} \to Y$ is a $[L]$-fibration.
\end{itemize}
\end{lem}
\begin{proof}
{\bf (a)} Let ${\mathcal F}$ be the set of all inclusions
$\{ f_{t} \colon A_{t} \hookrightarrow B_{t} ; t \in T\}$ such that $B_t$ is an $n$-$[L]$ disk and $A_{t}$ is a corresponding $(n-1)$-$[L]$-sphere, $n=0,1,\dots$, or
$A_{t}$ and $B_{t}$ are finite $[L]$-polyhedra such that both $A_t$ and $B_t$ are $AE([L])$--spaces.
For each $t \in T$ consider the set $S(t)$ which contains all
pairs of maps $(g,h)$ such that
the following
diagram
\medskip

\[
\begin{diagram}
\node{A_{t}} \arrow{e,t}{g} \arrow{s,l}{f_{t}}  \node{X} \arrow{s,r}{p}\\
\node{B_{t}} \arrow{e,t}{h}  \node{Y}
\end{diagram}
\]
\medskip

\noindent commutes. By gluing a copy of $B_{t}$ to $X$ along
$A_{t}$ for every commutative diagram of the above form we
obtain the Gluing Construction $G^{1}({\mathcal F},p)$ as in
the following pushout diagram
\medskip

\[
\begin{diagram}
\node{\coprod\left\{\coprod\{ A_{t} \colon (g,h) \in S(t)\}
\colon t \in {\mathcal T}\right\}} \arrow{e,t}{\widetilde{g}}
\arrow{s,l}{\coprod\{ f_{t}\colon t \in T\}}  \node{X} \arrow{s,r}{p}\\
\node{\coprod\left\{\coprod\{ B_{t} \colon (g,h) \in S(t)\}
\colon t \in {\mathcal T}\right\}} \arrow{e,t}{\widetilde{h}}
\node{G^{1}({\mathcal F},p)}
\end{diagram}
\]
\medskip

\noindent Let $i_{1} \colon X \to G^{1}({\mathcal F},p)$ denote
the natural embedding. By the universal property of colimits,
there exists a map $p_{1} \colon G^{1}({\mathcal F},p) \to Y$
such that $p_{1}i_{1} = p$. By repeating this process we obtain the following
infinite commutative diagram
\medskip

\[
\begin{CD}
X @>i_{1}>> G^{1}({\mathcal F},p) @>i_{2}>> G^{2}({\mathcal F},p) @>i_{3}>> \cdots @>i_{k}>> G^{k}({\mathcal F},p) @> >> \cdots \\
@V{p}VV @V{p_{1}}VV @V{p_{2}}VV @VV{}V @V{p_{k}}VV    \\
Y @>\operatorname{id}>> Y @>\operatorname{id}>> Y @>\operatorname{id}>> \cdots @>\operatorname{id}>> Y @> >> \cdots ,
\end{CD}
\]
\medskip

\noindent Next consider the colimit $G^{\infty}({\mathcal F},p)$ of the
upper row (i.e. the Infinite Gluing Construction) in the above
diagram. Obviously there are natural maps
$i_{\infty} \colon X \to G^{\infty}({\mathcal F},p)$ and
$p_{\infty} \colon G^{\infty}({\mathcal F},p) \to Y$ such that
$p = p_{\infty}i_{\infty}$.

It follows from the small object argument (see \cite[Proposition
7.17]{dwyspa95}) and the proof of \cite[Lemma 8.12]{dwyspa95})
that the map $p_{\infty} \colon G^{\infty}({\mathcal F},p) \to Y$
has the RLP with respect to each of the maps in the family
${\mathcal F}$. This, according to definition and by Lemma
\ref{L:acyc}, means that $p_{\infty}$ is an acyclic
$[L]$-fibration.

Since every map $f_{t} \in {\mathcal F}$ is either an inclusion of
finite $[L]$-polyhedra, which are $AE([L])$--spaces,
or an inclusion of $(n-1)^{th}$-$[L]$--sphere into
$n^{th}$-$[L]$--disk for some $n$, it follows from Definition
\ref{D:Lmodel}(ii)$_{L}$ and Lemma \ref{L:acyc} that $f_{t}$ has
the LLP with respect to acyclic $[L]$-fibrations. The colimit
universality property of the Infinite Gluing Construction
guarantees that the map $i_{\infty}$ also has the LLP with respect
to all acyclic $[L]$-fibrations and therefore is a
$[L]$-cofibration according to Definition
\ref{D:Lmodel}(iii)$_{L}$.

{\bf (b)} Let ${\mathcal F}$ be the set of all inclusions
$\{ f_{t} \colon A_{t} \hookrightarrow B_{t} ; t \in T\}$ of
finite $[L]$-polyhedra such that $A_t$ and $B_t$ are $AE([L])$--spaces.
For each $t \in T$ consider the set $S(t)$ which contains all
pairs of maps $(g,h)$ such that
the following
diagram
\[
\begin{diagram}
\node{A_{t}} \arrow{e,t}{g} \arrow{s,l}{f_{t}}  \node{X} \arrow{s,r}{q}\\
\node{B_{t}} \arrow{e,t}{h}  \node{Y}
\end{diagram}
\]

\noindent commutes. As in the part (b), we obtain
the Infinite Gluing Construction
$G^{\infty}({\mathcal F},q)$ and natural maps
$j_{\infty} \colon X \to G^{\infty}({\mathcal F},q)$ and
$q_{\infty} \colon G^{\infty}({\mathcal F},q) \to Y$ such that
$q = q_{\infty}j_{\infty}$.

It follows from the small object argument (see
\cite[Proposition 7.17]{dwyspa95}) and the proof of
\cite[Lemma 8.12]{dwyspa95}) that the map
$q_{\infty} \colon G^{\infty}({\mathcal F},q) \to Y$ has
the RLP with respect to each of the maps in the family
${\mathcal F}$. This, according to definition, means that
$q_{\infty}$ is a $[L]$-fibration.

Since every map $f_{t} \in {\mathcal F}$
is an inclusion of finite $[L]$-polyhedra which are $AE([L])$--spaces,
it follows that $f_{t}$ has the
LLP with respect to $[L]$-fibrations (Definition
\ref{D:Lmodel}(ii)$_{L}$). It is clear that $j_{\infty}$
is a weak $[L]$-homotopy equivalence. The colimit universality
property of the Infinite Gluing Construction implies that the
map $j_{\infty}$ also has the LLP with respect to $[L]$-fibrations.
\end{proof}

\begin{thm}\label{T:main}
Let $L$ be a finite polyhedron. Then the category
${\mathbf T}{\mathbf O}{\mathbf P}$ with weak equivalences,
fibrations and cofibrations, defined as in
Definition \ref{D:Lmodel}, is a model category.
\end{thm}
\begin{proof}
Axioms (MC1)--(MC2) are trivially satisfied.
Next, observe that the notions of $[L]$--fibrations and $[L]$--cofibrations are defined by lifting properties and hence the classes of $[L]$--fibrations and $[L]$--cofibrations are closed by retracts. Further, a retract of an isomorphism is an isomorphism.
These facts imply Axiom (MC3).
Lemma \ref{L:factor} shows that Axiom (MC5) is also satisfied.

In order to verify Axiom (MC4) we need only to check that acyclic
$[L]$--cofibrations have the LLP with respect to all
$[L]$--fibrations ($[L]$--cofibrations have the property, required
in (MC4), by definition). Consider a commutative diagram
\medskip

\[
\begin{diagram}
\node{A} \arrow{e,t}{f} \arrow{s,l}{i}  \node{X} \arrow{s,r}{p}\\
\node{B} \arrow{e,t}{g} \arrow{ne,t,..}{h} \node{Y}
\end{diagram}
\]
\medskip

\noindent of unbroken arrows, where $i$ is an acyclic $[L]$--cofibration and $p$ is a $[L]$--fibration. We need to show existence of $h$.
By Lemma \ref{L:factor} we have the following commutative diagram
\medskip

\[
\begin{diagram}
\node{A} \arrow{e,t}{j_{\infty}} \arrow{s,l}{i}  \node{A'} \arrow{s,r}{q_{\infty}}\\
\node{B} \arrow{e,t}{\operatorname{id}} \arrow{ne,t,..}{h} \node{B}
\end{diagram}
\]
\medskip

\noindent where $j_{\infty}$ is a weak $[L]$--homotopy equivalence having the LLP with respect to all $[L]$--fibrations and $q_{\infty}$ is an $[L]$--fibration.

Since $i$ and $j_{\infty}$ are weak $[L]$--homotopy equivalences it follows from (MC2) that $q_{\infty}$ is a weak $[L]$--homotopy equivalence. Therefore $q_{\infty}$ is an acyclic $[L]$--fibration. Since $i$ is a $[L]$--cofibration, there exists a lifting $h'\colon B\to A'$ of
$q_{\infty}$. Since $j_{\infty}$ has the LLP with respect to $[L]$--fibrations, there exists a lifting $h''\colon A'\to X$ of a composition $g\circ q_{\infty}\colon A'\to Y$.
We let $h=h''\circ h'$. It is easily seen that $h$ has the required property.
\end{proof}

Theorem \ref{T:main} and Corollary \ref{C:cofibrant} imply the following important observation.

\begin{cor}\label{C:whitehead}
A map between $[L]$--complexes is an $[L]$--homotopy equivalence if and only if it induces isomorphisms of all $[L]$--homotopy groups.
\end{cor}

In conclusion let us show that our notion of $[S^{n}]$--fibration differs from the usual notion of $n$--fibration. Consequently the model category structure generated by $L = S^{n}$ in theorem \ref{T:main} differs from the one described in Theorem C (note, however, that these two structures have identical weak equivalences).
It is unclear if the model category structure generated on ${\bf TOP}$ by $L = \{ \text{{\rm point}}\}$ coincides with the one described in Theorem A (although classes of weak equivalences are identical and the notion of $[L]$-homotopy coincides with the notion of usual homotopy in this case).

To this end let $n \geq 0$ and let $X$ and $Y$ be copies of the $n$-dimensional universal Menger compactum $\mu^{n}$. Consider the Dranishnikov's resolution $p \colon X \to Y$ constructed in \cite{DR1} (see also \cite[\S 4.2]{chibook}).
Observe that the map $p$, being polyhedrally $n$--soft, is an acyclic $n$--fibration. Note that $p$ is $(n-1)$-soft, but not $n$--soft.
The latter means that there exists an at most $n$-dimensional compactum $B_{0}$, a closed embedding $i_{0} \colon A_{0} \hookrightarrow B_{0}$ and two maps $\alpha \colon A_{0}\to X$ and $\beta \colon B_{0} \to Y$ such that there is no lifting of $\beta$, extending $\alpha$. In other words, the following commutative diagram of unbroken arrows

\[
\begin{diagram}
\node{A_{0}} \arrow{e,t}{\alpha} \arrow{s,l}{i}  \node{X} \arrow{s,r}{p}\\
\node{B_{0}} \arrow{e,t}{\beta} \arrow{ne,t,..}{\widetilde{\alpha}} \node{Y}
\end{diagram}
\]
\medskip

\noindent cannot be completed by the broken arrow.

Let $A$ and $B$ denote two additional copies of the Menger compactum $\mu^{n}$ and let $i \colon A \to B$ be a $Z$--embedding.
Consider embeddings $j_{0}\colon A_{0} \to A$ and $j \colon B_{0} \to B$ such that $j\circ i_{0} = i\circ j_{0}$ and $j(B_{0}\setminus i_{0}(A_{0})) \subseteq B \setminus i(A)$ (embeddings with the indicated properties exist because $i$ is a $Z$-embedding).
Since $X$ is an ${\rm AE}(n)$--compactum and $\dim A = n$, it follows that there exists a map $g\colon A\to X$ such that $\alpha = g \circ i_{0}$.
Define the map $h \colon i(A) \cup j(B_{0}) \to Y$ by letting

\[
h(b) =
\begin{cases}
p(g(i^{-1}(b))), & \text{if}\; b \in i(A),\\
\beta (j^{-1}(b)), & \text{if}\; b \in j(B_{0}).
\end{cases}
\]

\noindent Since $Y\in {\rm AE}(n)$ and $\dim B = n$, there exists a map
$f \colon B \to Y$ such that $f|\left(i(A) \cup j(B_{0})\right) =h$. Note that $p\circ g = f \circ i$.

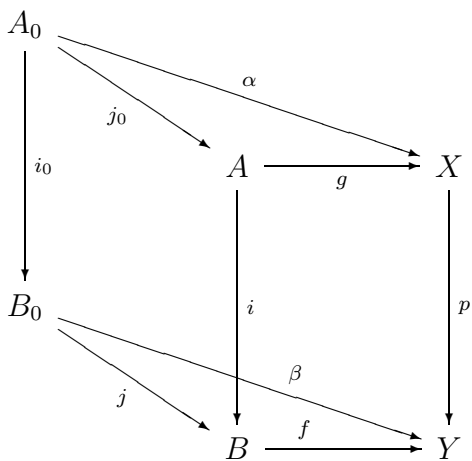
\[
\dgARROWLENGTH=2.0\dgARROWLENGTH
\dgARROWPARTS=8
\begin{diagram}
\node{A_{0}} \arrow{se,r}{j_{0}}  \arrow{ese,l}{\alpha} \arrow[2]{s,b}{i_{0}}\\
\node[2]{A}  \arrow[2]{s,r}{i} \arrow{e,r}{g} \node{X} \arrow[2]{s,r}{p}\\
\node{B_{0}} \arrow{ese,l,5}{\beta} \arrow{se,b}{j}\\
\node[2]{B} \arrow{e,t,2}{f} \node{Y}
\end{diagram}
\]

\medskip

Finally observe that the map $f$
does not have a lifting extending $g$, since this would imply existence of lifting of $\beta$ extending $\alpha$. Hence $p$ is not an $[S^{n}]$--fibration.



\end{document}